\documentclass[12pt,a4paper]{article}
\usepackage[cp850]{inputenc}
\usepackage{latexsym, amsfonts, amssymb}
\usepackage[english]{babel}
\usepackage{latexsym}
\usepackage{amscd}
\newcommand{\beq}{\begin{equation}}
\newcommand{\eeq}{\end{equation}}
\newcommand{\beas}{\begin{eqnarray*}}
\newcommand{\eeas}{\end{eqnarray*}}
\newcommand{\qq}{\mathbb{Q}}
\renewcommand{\P}{{\mathcal P}}
\newcommand{\D}{{\cal D}}
\renewcommand{\O}{{\cal O}}
\newcommand{\X}{{\cal X}}
\newcommand{\Y}{{\cal Y}}
\renewcommand{\S}{{\mathcal S}}
\newcommand{\gen}[1]{\mbox{$\langle #1 \rangle$}}
\newcommand{\sgn}{\textrm{sgn}}

\newtheorem{theorem}{Theorem}
\newtheorem{lemma}[theorem]{Lemma}
\newtheorem{corollary}[theorem]{Corollary}
\newtheorem{proposition}[theorem]{Proposition}

\begin{document}
\begin{center}
{\Large\bf Properties of some character tables}\\[5pt]
{\Large\bf related to the symmetric groups}\\[.1in]
{\large Christine Bessenrodt, J\o rn B. Olsson$^1$, Richard
  P. Stanley$^2$}\\[1ex]
\end{center}

\medskip

\begin{abstract} 
We determine invariants like the Smith normal form
and the determinant for certain integral matrices which arise from the
character tables of the symmetric groups $S_n$ and their double
covers. 
In particular, we give a simple computation, based on the
theory of Hall-Littlewood symmetric functions, of the determinant of
the regular character table ${\cal X}_{RC}$ of $S_n$ with respect to an
integer $r\geq 2$. This result had earlier been proved by Olsson in a
longer and more indirect manner. 
As a consequence, we obtain a  new proof of the Mathas' Conjecture 
on the determinant of the Cartan matrix of the Iwahori-Hecke algebra. 
When $r$ is prime we determine the
Smith normal form of ${\cal X}_{RC}$. Taking $r$ large yields the Smith
normal form of the full character table of $S_n$. Analogous results
are then given for spin characters. 
\end{abstract}

\noindent Christine Bessenrodt \\
Institut f\"ur Mathematik\\
Universit\"at Hannover\\
D-30167 Hannover, Germany\\
bessen@math.uni-hannover.de\\
\smallskip

\noindent J\o rn B. Olsson\\
Matematisk Afdeling\\
University of Copenhagen \\
Copenhagen, Denmark\\
olsson@math.ku.dk\\

\noindent Richard P. Stanley\\
Department of Mathematics 2-375\\
M.I.T.\\
Cambridge, MA 02139, USA\\
rstan@math.mit.edu\\

\noindent$^1$Partially supported by The Danish National Research Council.

\noindent $^2$Partially supported by NSF grant \#DMS-9988459.

\pagebreak

\section{Introduction}\label{Intro}

In this paper we determine invariants like the Smith normal form or
the determinant for certain integral matrices which come from the
character tables of the finite symmetric groups $S_n$ and their double
covers $\hat{S}_n.$ The matrices in question are the so-called regular
and singular character tables of $S_n$ and the reduced spin character
table of $\hat{S}_n.$

In section~\ref{det} we calculate the determinants of the $r$-regular and
$r$-singular character tables of $S_n$ for arbitrary integers $r \ge
2,$ using symmetric functions and some bijections involving regular
partitions. The knowledge of these determinants is equivalent to the
knowledge of the determinants of certain ``generalized Cartan
matrices'' of $S_n$ as considered in \cite{KOR}. In particular we
obtain a new proof of a conjecture of Mathas about the Cartan
matrix of an Iwahori-Hecke algebra of $S_n$ at a primitive $r$th root
of unity which is simpler than the original proof given by Brundan and
Kleshchev in \cite{BruKle}.
In section~\ref{Smith} we determine the Smith normal form of the
regular character table in the case where $r$ is a prime. As a special
case the Smith normal form of the character table of $S_n$ may be
calculated. We also determine the Smith normal form of the reduced
spin character table for $\hat{S}_n.$ The paper also presents some
open questions.

\bigskip

\section{The determinant of the regular part of the character table of
$S_n$}\label{det}

We fix positive integers $n,r$, where $r \ge 2.$

If $\mu=(\mu_1,\mu_2, \ldots )$ is a partition of $n$ we write $\mu \in
\P$ and denote by $\ell(\mu)$ the number of (non-zero) parts of $\mu$.
We let $z_{\mu}$ denote the order of the centralizer of an
element of (conjugacy) type $\mu$ in $S_n$. Suppose
$\mu=(1^{m_1(\mu)}, 2^{m_2(\mu)},\ldots)$, is written in exponential
notation.  Then we may factor $z_\mu= a_{\mu}b_{\mu}$, where
$$a_{\mu}=\prod_{i \ge 1}i^{m_i(\mu)}, \hspace{.1 in}
b_{\mu}=\prod_{i \ge 1}m_i(\mu)!$$

Whenever $\mathcal{Q} \subseteq \mathcal{P}$
we define
$$a_{\mathcal{Q}}=\prod_{\mu \in \mathcal{Q}}a_{\mu}, \hspace{.1 in}
b_{\mathcal{Q}}=\prod_{\mu \in \mathcal{Q}}b_{\mu}.$$

Let $\mu \in \P.$
We write $\mu \in  R$ and call $\mu$  {\it regular} if
$m_i(\mu) \le r-1$ for all $i \ge 1.$
We write $\mu \in  C$ and call $\mu$  {\it class regular} if
$m_i(\mu)=0,$ whenever $r \mid i.$

We are particularly interested in the integers $a_C$ and  $b_C$.
By \cite[Theorem 4]{OlsReg} there is a connection between $a_C$ and
$b_C$  given by
\beq b_C=r^{d_C} \,a_C,
\label{ab-rel}\eeq
where the class regular defect number $d_C$ is defined by
$$d_C=\sum_{\mu \in C}d(\mu),
~~~d(\mu)=\sum_{i,k \ge 1}\left\lfloor \frac{m_i(\mu)}{r^k}\right\rfloor.$$
Here $\lfloor \cdot \rfloor$ is the floor function, i.e., 
 $\lfloor x \rfloor$ denotes the integral part of $x$. 
Note that for $r>n$ we have $R=C=\P$ and then $d_{\P}=0$ and thus
$a_{\P}=b_{\P}$.

\medskip
Let $\mathcal{X}_{RC}$ denote the regular character table of $S_n$
with respect to $r$. It is a submatrix of the character table
$\mathcal{X}$ of $S_n.$ The subscript $RC$ indicates that the rows of
$\mathcal{X}_{RC}$ are indexed by the set $R$ of regular partitions of
$n$, and the columns by the set $C$ of class regular partitions of
$n$. We want to present a proof of the following result:

\medskip

\begin{theorem}\label{det-regchartable}
We have
  $$ |\det \left(\mathcal{X}_{RC}\right)| = a_C. $$
\end{theorem}

This result was first proved in \cite{OlsReg}, but the proof relied
on results of \cite{KOR} for which the work of Donkin \cite{Donkin} and
Brundan and Kleshchev \cite{BruKle} was used in a crucial way. 
Our proof of Theorem~\ref{det-regchartable} 
does not use  \cite{Donkin} or \cite{BruKle}; 
it is direct and thus much shorter. 

\smallskip 

In \cite{KOR}, an $r$-analogue of the modular
representation theory for $S_n$ was developed systematically, and 
in particular, an $r$-analogue of the Cartan matrix
for the symmetric groups (and the corresponding $r$-blocks)  
was introduced.

 In \cite{BesOls} the explicit value of this latter 
determinant was conjectured to be $r^{d_C}$ in the notation above; 
this was proved in \cite[Proposition 6.11]{KOR} 
using \cite{Donkin} and \cite{BruKle}.    
This result is now a consequence of our theorem:

\begin{corollary}\label{KORdet} 
Let $\mathcal{C}$ be the $r$-analogue of the Cartan matrix of $S_n$ 
as defined in \cite{KOR}. 
Then we have
 $$ \det \left(\mathcal{C}\right) = r^{d_C} \: . $$
\end{corollary}

\textbf{Proof.} As is shown in \cite{OlsReg} there is a simple 
equation connecting the determinants of $\mathcal{C}$ and 
$\mathcal{X}_{RC}$,  namely
$$ \det(\mathcal{X}_{RC})^2 \det(\mathcal{C}) =a_Cb_C\:.$$
Thus in view of equation (\ref{ab-rel}) 
Theorem~\ref{det-regchartable} implies the 
Corollary. 
$\ \Box$

\medskip

Mathas conjectured that the determinant of the Cartan matrix of
an Iwahori-Hecke algebra of $S_n$ at a primitive $r$th root of unity
should be a power of $r$; 
via \cite{Donkin},  the conjecture in \cite{BesOls} mentioned above  
predicted the explicit value of this determinant,  
thus providing a strengthening of Mathas' conjecture.
Mathas' conjecture  was proved by Brundan and Kleshchev
\cite{BruKle};  in fact, they 
also gave an explicit formula for this determinant 
for blocks of the Hecke algebra. 
We can now provide an alternative proof of  
these conjectures. 

\begin{corollary}\label{Mathas}
The strengthened Mathas' conjecture is true.
\end{corollary}

\textbf{Proof.} Donkin \cite{Donkin} has shown that the Cartan matrix 
for the Hecke algebra has the same determinant 
as  the Cartan matrix $\mathcal{C}$ considered in Corollary~\ref{KORdet}.
$\ \Box$
 
\medskip

Based on this and the results on $r$-blocks in \cite{KOR},  
the results in \cite{BesOls} then also give the determinants 
of Cartan matrices of $r$-blocks of $S_n$ explicitly, 
without the use of \cite{BruKle}.
 
Let us finally mention that in \cite[Section 6]{KOR}
there is an explicit conjecture about the Smith normal form  of
$\mathcal{C}$. 
In the case where $r$ is a prime, this is known to be true by the
general theory of R.\ Brauer. 
One may also ask about the Smith normal
form of $\mathcal{X}_{RC}$; 
we answer this question in this article in the prime case.

\medskip

We now proceed to describe the proof of Theorem~\ref{det-regchartable}.
It is obtained by
combining Theorems~\ref{det^2-regchartable} and~\ref{B/A}  below.
Theorem~\ref{det^2-regchartable}
 evaluates $\det \left(\mathcal{X}_{RC}\right)^2$ using
symmetric functions as an expression involving a primitive $r$th
root of unity. Theorem~\ref{B/A} shows that this expression equals ${a_C}^2$.
It is based on general bijections involving regular
partitions.

Define
   $$ z_\mu(t) = z_\mu \prod_j(1-t^{\mu_j})^{-1}=
z_\mu \prod_i (1-t^i)^{-m_i(\mu)} $$
where the product ranges over all $j$ for which $\mu_j>0$, and
  $$ b_\lambda(t) = \prod_i (1-t)(1-t^2)\cdots (1-t^{m_i(\lambda)}).
  $$
Let $\omega=e^{2\pi i/r}$, a primitive $r$th root of unity.

We use notation from the theory of symmetric functions from
\cite{macd} or \cite{ec2}. In particular, $m_\lambda$, $s_\lambda$,
and $p_\lambda$ denote the monomial, Schur, and power sum symmetric
functions, respectively, indexed by the partition $\lambda$.

\medskip

\begin{theorem}\label{det^2-regchartable}
We have
  $$ \det \left(\mathcal{X}_{RC}\right)^2 = \prod_{\mu\in
    C}z_\mu(\omega)\cdot
    \prod_{\lambda\in R} b_\lambda(\omega). $$
\end{theorem}

\textbf{Proof.}
Let $Q_\lambda(x;t)$ denote a Hall-Littlewood
symmetric function as in \cite[p.\ 210]{macd}. It is immediate from
the definition of $Q_\lambda(x;t)$ that $Q_\lambda(x;\omega)=0$ unless
$\lambda\in R$. Moreover (see \cite[Exam.~III.7.7,
p.~249]{macd}) when $Q_\lambda(x;\omega)$ is expanded in terms of
power sums $p_\mu$, only class regular $\mu$ appear. Thus
\cite[(7.5), p.~247]{macd} for $\lambda\in R$ we have
  $$ Q_\lambda(x;\omega) = \sum_{\mu\in C}
  z_\mu(\omega)^{-1}X_\mu^\lambda(\omega)p_\mu(x), $$
where $X_\mu^\lambda(t)$ is a Green's polynomial.

Hence by \cite[(7.4)]{macd} the matrix $X(\omega)_{RC} =
(X_\mu^\lambda(\omega))$,
where $\lambda\in R$ and $\mu\in C$, satisfies
  \beq \det(X(\omega)_{RC})^2 = \prod_{\mu\in C}
      z_\mu(\omega) \prod_{\lambda \in R}b_\lambda(\omega).
     \label{xdet} \eeq

Now consider the symmetric function $S_\lambda(x;t)$ as defined in
\cite[(4.5), p.~224]{macd}. It follows from the formula
$S_\lambda(x;t)=s_\lambda(\xi)$ in \cite[top of p.~225]{macd} that
  $$ S_\lambda(x;t) = s_\lambda(p_j\rightarrow (1-t^j)p_j), $$
i.e., expand $s_\lambda(x)$ as a polynomial in the $p_j$'s and
substitute $(1-t^j)p_j$ for $p_j$. Since
  $$ s_\lambda = \sum_\mu z_\mu^{-1}\chi^\lambda(\mu) p_\mu, $$
we have
  $$ S_\lambda(x;\omega) = \sum_{\mu\in C} z_\mu(\omega)^{-1}
      \chi^\lambda(\mu)p_\mu. $$
The $S_\lambda(x;\omega)$'s thus lie in the space $A_{(r)}$ spanned
over $\qq(\omega)$ by the $p_\mu$'s where $\mu\in C$. 
Since the
$Q_\mu(x;\omega)$'s for regular $\mu$ span $A_{(r)}$ by
\cite[Exam.~III.7.7, p.~249]{macd}, the same is true of the
$S_\lambda(s;\omega)$'s. 
Moreover, the
transition matrix $M(S,Q)_{RR}$ between the $Q_\lambda(x;t)$'s and
$S_\lambda(x;t)$'s is lower unitriangular by \cite[top of
p.~239]{macd} and \cite[p.~241]{macd}. 
Hence
  \beq \det M(S,Q)_{RR}=1. \label{msq} \eeq

Let $M(S,p)_{RC}$ denote the transition matrix from the $p_\mu$'s to
$S_\lambda$'s for $\mu\in C$ and $\lambda\in R$. Let $Z(t)_{CC}$
denote the diagonal matrix with entries $z_\lambda(t)$, $\lambda\in
C$. By the discussion above we have
  \beas \mathcal{X}_{RC} & = & M(S,p)_{RC}\,Z(\omega)_{CC}
  \ \ \mbox{(by the relevant definitions)}\\ & = &
     M(S,Q)_{RR}\,M(Q,p)_{RC}\,Z(\omega)_{CC}\\ & = &
    M(S,Q)_{RR}\,X(\omega)_{RC}\,Z(\omega)^{-1}_{CC}\,Z(\omega)_{CC}\\
    & = & M(S,Q)_{RR}\,X(\omega)_{RC}. \eeas
Taking determinants and using (\ref{xdet}) and (\ref{msq}) completes
the proof. $\ \Box$

\bigskip

Define $$A_C(\omega)=\prod_{\mu \in C}\prod_{i}(1-\omega^i)^{-m_i(\mu)}$$
$$B_R(\omega)=\prod_{\lambda \in R}b_{\lambda}(\omega)^{-1}=
\prod_{\lambda \in R}
\left(\prod_i (1-\omega)(1-\omega^2)\cdots (1-\omega^{m_i(\lambda)})\right)^{-1},$$
so that by Theorem~\ref{det^2-regchartable}
 $$ \det \left(\mathcal{X}_{RC}\right)^2 = a_Cb_C A_C(\omega)B_R(\omega)^{-1}. $$

In order to complete the proof of Theorem~\ref{det-regchartable} we thus just need to show:
$$ \frac{B_R(\omega)}{A_C(\omega)}=\frac{b_C}{a_C} \: .$$
As ~~$\displaystyle\frac{b_C}{a_C} = r^{d_C}$ ~~this is equivalent to showing

\begin{theorem}\label{B/A}
We have
$$ \frac{B_R(\omega)}{A_C(\omega)}=r^{d_C}.$$
\end{theorem}

Clearly the factors $1-\omega^j$ occurring on the left hand side in Theorem~\ref{B/A}
depend only on the residue of $j$ modulo $r.$ Thus
$$A_C(\omega)^{-1}=\prod_{s=1}^{r-1}(1-\omega^s)^{\alpha_C^{(s)}}, ~~~
B_R(\omega)^{-1}=\prod_{s=1}^{r-1}(1-\omega^s)^{\beta_R^{(s)}}, $$
where
$$ \alpha_C^{(s)}=\sum_{\mu \in C}
\sum_{\{i | i \equiv s (\mathrm{mod}~ r) \}} m_i(\mu)$$
 $$ \beta_R^{(s)}=\sum_{\rho \in R} |\{i|m_i(\rho) \ge s \}|. $$

We use the bijections $\kappa^{(s)}$ defined in Proposition~\ref{Prop7} 
below to show the following:

\medskip

\begin{proposition}\label{Prop4}
For all $s \in \{1, \ldots , r-1\}$ we have
$$\alpha_C^{(s)}=\beta_R^{(s)}+d_C.$$
\end{proposition}

This shows then that
$$ \frac{B_R(\omega)}{A_C(\omega)}= \left(\prod_{s=1}^{r-1}(1-\omega^s)\right)^{d_C}.$$
Then Theorem~\ref{B/A}  follows from the fact that
$$\prod_{s=1}^{r-1}(1-\omega^s)=r.$$
(Simply substitute $x=1$ in the identity
$1+x+\cdots+x^{r-1}=\prod_{s=1}^{r-1}(x-\omega^s)$.)

\bigskip

Let $m \in \mathbb{N}$.
We write $m$ in its $r$-adic decomposition as
$m=\sum_{j\geq 0}m_jr^j$, i.e., with
 $m_j\in \{ 0, \ldots, r-1\}$ for all~$j$.
For  $m \neq 0$, we can write
$m=\sum_{j\geq k}m_jr^j$, with $m_k\neq 0$.
In the power series convention, $k(m)=k$ is the degree of~$m$
and $\ell(m)=m_k$ its leading coefficient.
We also set $h(m)=\sum_{j\geq k+1}m_jr^j = r^{k+1}q(m)$
for the higher terms of~$m$.
Thus
$$m=\ell(m)  r^{k(m)} + q(m) r^{k(m)+1} \: .$$
For a given $a$, we define
$$h_a(m)=\sum_{j\geq a}m_jr^j = q_a(m)  r^a ~,~~
q_a(m)=\left\lfloor \frac{m}{r^a} \right\rfloor \:.$$

We call $e \in \{1, \ldots , m\}$ a {\it non-defect number for~$m$},
if $h(e)=h_{k(e)+1}(m)$, otherwise $e$ is a {\it defect number for~$m$}
(and then $h(e)< h_{k(e)+1}(m)$, and hence $q(e)< q_{k(e)+1}(m)$).
Thus the non-defect numbers for $m$ are of the form
$$e=e_a r^a + h_{a+1}(m) ~,~~ e_a \in \{1, \ldots , m_a\} \: , $$
and thus there are $\sum_{j\geq 0}m_j$ such numbers.
The defect numbers for $m$ are of the form
$$e=e_ar^a+qr^{a+1}
~,~~ e_a\in\{1, \ldots,r-1\}~,~ q\in \{0, \ldots , q_{a+1}(m)-1\}
\:.$$
Their parameters $(a,q)$ thus belong to the set
$$\mathcal{D}(m)=\{(a,q) \mid a \ge 0,~ 0 \le q <
q_{a+1}(m) \}, $$
which is of cardinality
$$ d(m)=\sum_{a \ge 1} \left\lfloor \frac{m}{r^a}\right\rfloor \:,$$
called the {\it defect} of $m$.
For each $s \in \{1, \ldots , r-1\}$ there are exactly
$d(m)$ defect numbers for $m$ with leading coefficient $s$,
namely $e=s r^a + q r^{a+1}$, where $(a,q)\in \mathcal{D}(m)$.
Thus clearly we have $(r-1)d(m)$ defect numbers for $m$ and
$$ m=(r-1)d(m)+\sum_{j \ge 0} m_j \: .$$

\bigskip

For $\mu \in \P$, its defect (as defined at the beginning of this section)
is then
$$ d(\mu)=\sum_{i \ge 1} d(m_i(\mu)).$$
For $s \in \{1, \ldots , r-1\}$ set
$$ \mathcal{D}^{(s)}(\mu)=\{(i,a,q) \mid \ell(i)=s, (a,q) \in
\mathcal{D}(m_i(\mu)) \}  $$
and
$$\mathcal{D}(\mu)=\bigcup_{s=1}^{r-1} \mathcal{D}^{(s)}(\mu) \: .$$
We have that
$$  d^{(s)}(\mu)=|\mathcal{D}^{(s)}(\mu)|=
\sum_{\{i \ge 1, \ell(i)=s\}} d(m_i(\mu) )$$
and
$$ d(\mu)= \sum_{s=1}^{r-1}d^{(s)}(\mu) = |\mathcal{D}(\mu)| \: .$$

Consider nonzero residues $s,t$ modulo $r$, let $\mu=(i^{m_i(\mu)})$
and define
$$\mathcal{T}^{(st)}(\mu)= \{(i,j) \mid 1 \leq i,~1 \leq j \leq m_i(\mu),
\ell(i)=s, \ell(j)=t \}.$$

\medskip
Glaisher \cite{Glaisher} defined
a bijection between the sets $C$ and $R$ of class
regular and regular partitions of $n.$ Glaisher's map $G$ is defined
as follows. Suppose that $\mu=(i^{m_i(\mu)}) \in C.$
Consider the  $r$-adic expansion of each multiplicity $m_i(\mu):$
$$m_i(\mu)=\sum_{j \ge 0} m_{ij}(\mu)r^j$$
where for all relevant $i,j$ we have $m_{ij}(\mu) \in \{0, \ldots ,r-1\}$.
Then $G(\mu)=\rho$
where for all $i,j, r \nmid i$ we have $m_{ir^j}(\rho)=m_{ij}(\mu).$

\smallskip

We show

\begin{proposition}\label{Prop5}
If $\mu \in C$ then
$|\mathcal{T}^{(st)}(\mu)|=|\mathcal{T}^{(st)}(G(\mu))|+d^{(s)}(\mu).$
\end{proposition}

\smallskip

\textbf{Proof.} We establish a bijection $\delta^{(st)}(\mu)$ between
$\mathcal{T}^{(st)}(\mu)$ and the disjoint union
$ \mathcal{T}^{(st)}(G(\mu)) \cup \mathcal{D}^{(s)}(\mu).$
If $(i,j) \in \mathcal{T}^{(st)}(\mu)$ and $(k(j), q(j))=(a,q)$,
we have two possibilities

(i) $j$ is a defect number for $m_i(\mu)$.
Then we map $(i,j)$ onto
$(i,a,q) \in \mathcal{D}^{(s)}(\mu).$

(ii) We have
$j= t r^a + h_{a+1}(m_i(\mu))$ where $1 \le t \le
m_{ia}(\mu).$ Then we map $(i,j)$ onto $(r^ai,t) \in \mathcal{T}^{(st)}(G(\mu)). $

This establishes the desired bijection. $\ \Box$

\medskip

Consider nonzero residues $s,t$ modulo $r$, and define
$$\mathcal{T}^{(st)}_C= \{(\mu,i,j)| \mu \in C, (i,j) \in
\mathcal{T}^{(st)}(\mu)  \}$$
$$\mathcal{T}^{(st)}_R= \{(\rho,i,j)| \rho \in R, (i,j) \in
\mathcal{T}^{(st)}(\rho) \}$$
$$\mathcal{D}^{(s)}= \{(\mu,i,a,q)| \mu \in C,
(i,a,q)  \in \mathcal{D}^{(s)}(\mu)\}.$$

\medskip

Clearly the bijections $\delta^{(st)}(\mu),~ \mu \in C,$
above induce a bijection
$$\delta^{(st)}: \mathcal{T}^{(st)}_C
\longleftrightarrow \mathcal{T}^{(st)}_R \cup \mathcal{D}^{(s)}.$$

Putting the bijections $\delta^{(ts)},~t=1, \ldots, r-1$ together we
obtain a bijection

$$\delta^{(s)}: \bigcup_{t=1}^{r-1} \mathcal{T}^{(ts)}_C
\longleftrightarrow \bigcup_{t=1}^{r-1} \mathcal{T}^{(ts)}_R
\cup \mathcal{C},$$

where
$$\mathcal{C}=\bigcup_{t=1}^{r-1}\mathcal{D}^{(t)}.$$

\medskip

In \cite[proof of Theorem~4]{OlsReg}, an involution $\iota$ was
defined on the set
$$\mathcal{T}_C= \{(\mu,i,j)| \mu \in C, i,j \ge 1,
m_i(\mu) \ge j \}.$$

From the definition of $\iota$ it follows that it maps the subset
$\mathcal{T}^{(st)}_C$ of $\mathcal{T}_C$ into
$\mathcal{T}^{(ts)}_C. $ Thus we
conclude

\medskip

\begin{lemma}
For all
  $s \in \{1, \ldots , r-1\}$
there is a bijection
$$\iota^{(s)}: \bigcup_{t=1}^{r-1}\mathcal{T}^{(st)}_C
\longleftrightarrow  \bigcup_{t=1}^{r-1}\mathcal{T}^{(ts)}_C \;  .$$
\end{lemma}

\medskip

Composing the bijections $\iota^{(s)}$ and $\delta^{(s)}$ we see

\medskip

\begin{proposition}\label{Prop7}
For all   $s \in \{1, \ldots , r-1\}$
there is a bijection
$$\kappa^{(s)}:  \bigcup_{t=1}^{r-1}\mathcal{T}^{(st)}_C
\longleftrightarrow \bigcup_{t=1}^{r-1} \mathcal{T}^{(ts)}_R \cup
\mathcal{C}.  $$
\end{proposition}

{\bf Proof of Proposition~\ref{Prop4}.}
Just consider the cardinalities of the sets occurring in
Proposition~\ref{Prop7}. 
$$ \left|\bigcup_{t=1}^{r-1}\mathcal{T}^{(st)}_C\right| =
\sum_{\mu \in C} \sum_{\{i | \ell(i)=s \}} m_i(\mu) = \alpha_C^{(s)}.$$
The latter equality holds because a class regular partition contains no
parts divisible by $r.$ Thus if $m_i(\mu) \ne 0$ then $\ell(i)=s$
if and only if  $i \equiv s (\mathrm{mod}~ r) .$
$$\left|\bigcup_{t=1}^{r-1} \mathcal{T}^{(ts)}_R\right|=
\sum_{\rho \in R} |\{i|m_i(\rho) \ge s \}|= \beta_R^{(s)}. $$
This is because parts in regular partitions have multiplicities $<r.$
Finally 
$$|\mathcal{C}|=\sum_{t=1}^{r-1}d^{(t)}=\sum_{\mu \in C} d(\mu)=d_C.
\quad\Box $$

\medskip

\textbf{Remark.} There is of course also a {\it singular} character
table for $S_n,$  which we denote $\mathcal{X}_{R^{'}C^{'}}.$
It is also a submatrix of the character table $\mathcal{X}$ of $S_n.$
The subscript $R^{'}C^{'}$ indicates that the rows of
$\mathcal{X}_{R^{'}C^{'}}$ are
indexed by the set
$R^{'}$ of singular (i.e. nonregular) partitions of $n$,  and the
columns by the set $C^{'}$ of class singular (i.e. non-class regular)
partitions of $n$. For this we have
\begin{equation} |\det \left(\mathcal{X}_{R^{'}C^{'}}\right)| =
  b_{C^{'}}. \label{eq:rpcp} \end{equation}
There are different ways of proving this. In \cite{OlsReg} there is a
proof based on Theorem \ref{det-regchartable} and a result in
\cite{KOR}.

Another way of proving (\ref{eq:rpcp}) is \emph{via} an identity of
Jacobi \cite[p.\ 21]{gantmacher}. Namely, suppose that $A$ is an
invertible $n\times n$ matrix, and write $A$ and $A^{-1}$ in the block
form
  $$ A=\left[ \begin{array}{cc} B & C\\ D & E \end{array} \right],
    \quad A^{-1}=\left[ \begin{array}{cc} B' & C'\\ D' & E'
    \end{array} \right], $$
where $B$ and $B'$ are $k\times k$ matrices. Then
  $$ \det E' = \frac{\det B}{\det A}. $$
By the orthogonality of characters we have
   $$ {\cal X}^{-1} =  {\cal X}^t \Delta(z_\mu^{-1}), $$
where $\Delta(z_{\mu}^{-1})$ is the diagonal
matrix with the $z_{\mu}^{-1},~\mu \in \mathcal{P}$, on the diagonal.
Equation (\ref{eq:rpcp}) follows immediately from this observation and
Theorem~\ref{det-regchartable}.

\medskip

\textbf{Remark.} If we keep $r$ fixed and let $n$ vary, then the
result of Proposition~\ref{Prop4} may also be proved by calculating the
generating functions for $\alpha_C^{(s)}, \beta_R^{(s)}$ and $d_C.$
Indeed, if $P(q)$ is the generating function for the number of
partitions of $n,$ then $P_r(q)=\frac{P(q)}{P(q^r)}$ is the generating
function for the number of regular partitions of $n.$
We may then express the generating functions for  $\alpha_C^{(s)},
\beta_R^{(s)}$ and $d_C$ respectively by

$$A^{(s)}(q)=P_{r}(q)\sum_{i \ge 0} \frac{q^{ir+s}}{1-q^{ir +s}}$$
$$B^{(s)}(q)=P_{r}(q)\sum_{j \ge 1} \frac{q^{js}-q^{jr}}{1-q^{jr}}$$
$$D(q)=P_{r}(q)\sum_{j \ge 1} \frac{q^{jr }}{1-q^{jr}}.$$

We omit the details. From this Proposition~\ref{Prop4} may be deduced easily.

\section{Smith normal forms of character tables related to $S_n$}\label{Smith}

For a partition $\lambda$ of $n$, we denote by $\xi^{\lambda}$ 
the permutation character of $S_n$ obtained by inducing the trivial 
character of the Young subgroup $S_{\lambda}$ up to $S_n$. 
First we explicitly describe the values of these permutation characters  
(this is included here as we have not been able 
to find a reference for it).

\begin{proposition}\label{permcharvalues}
Let $\lambda, \mu \in \P$, $k=\ell(\lambda)$, $\ell = \ell(\mu)$. Then the
value $\xi^{\lambda}(\mu )$ of the permutation character $\xi^{\lambda}$ on the conjugacy
class of cycle type $\mu$ equals the number of ordered set partitions
$(B_1, \ldots , B_k)$ of $\{1, \ldots , \ell\}$ such that
$$\lambda_j = \sum_{i \in B_j} \mu_i \quad \mbox{for } \:  j \in \{1, \ldots , k\} \:.$$
\end{proposition}

\textbf{Proof.}
Let $\sigma_{\mu}$ be a permutation of cycle type~$\mu$.  Then (see
\cite{JK}) $\xi^{\lambda}(\mu )$ is the number of $\lambda$-tabloids
fixed by~$\sigma_{\mu}$.  Now clearly, a $\lambda$-tabloid is fixed
by~$\sigma_{\mu}$ if and only if its rows are unions of complete
cycles of $\sigma_{\mu}$. Thus such a decomposition of rows
corresponds to an ordered set partition $(B_1, \ldots , B_k)$ of the
cycles of $\mu$ with the sum conditions in the statement of the
Proposition.  $\ \Box$

\medskip

\textbf{Remark.}
One may also use a symmetric function argument for
computing the values  $R_{\lambda\mu}= \xi^{\mu}(\lambda )$.
The complete homogeneous symmetric function $h_\lambda$
is the (Frobenius) characteristic of the character
$\xi^\lambda$ (see \cite[Cor. 7.18.3]{ec2}), so
$h_\lambda = \sum_\mu z_\mu^{-1} R_{\lambda\mu} p_\mu$.
As the $h_\lambda$ and $m_\mu$ are dual bases,
as well as the $p_\lambda$ and $z_\mu^{-1} p_\mu$,
it then follows that $p_\lambda = \sum_\mu R_{\lambda\mu} m_\mu$.
Using \cite[Prop. 7.7.1]{ec2} then also  gives the formula
in Proposition~\ref{permcharvalues}.

\begin{corollary}\label{permchartable}
Let $\lambda , \mu \in \P$. Then we have
\begin{enumerate}
\item[{(i)}]
$\xi^{\lambda}(\mu ) = 0$ unless $\lambda \geq \mu$ (dominance order).
\item[{(ii)}]
$\xi^{\lambda}(\lambda ) = b_{\lambda}= \prod_i m_i(\lambda)!$.
\item[{(iii)}]
$\xi^{\lambda}(\lambda ) \mid \xi^{\lambda}(\mu )$.
\end{enumerate}
\end{corollary}

\textbf{Proof.}
Using the remark above, parts (i) and (ii) follow immediately
by~\cite[Cor.~7.7.2]{ec2} (or one may also prove it directly using
Proposition~\ref{permcharvalues}).\\
For (iii), we use the combinatorial description
given in Proposition~\ref{permcharvalues}. With notation as before,
let $(B_1,\dots,B_k)$ be an ordered partition
of the set $\{1, \ldots ,\ell\}$ contributing to
  $\xi^{\lambda}(\mu )$, i.e., satisfying the sum conditions.
Now any permutation of $\{1, \ldots ,k \}$ which interchanges only
parts of $\lambda$ of equal size leads to a permutation of the entries
of $(B_1,\dots,B_k)$ such that the corresponding ordered partition
still satisfies the sum conditions.  Hence $\xi^{\lambda}(\mu)$ is
divisible by $\prod_i m_i(\lambda)! = b_{\lambda}$ and thus by
$\xi^{\lambda}(\lambda)$.  $\ \Box$

\bigskip

We can now determine the Smith normal form for the regular character
table of $S_n$ in the case where $r=p$ is prime.

For an integer matrix $A$ we denote by $\S(A)$ its Smith normal form.
If $p$ is a prime, we write $A_{p'}$ for the matrix obtained by
taking only the $p'$-parts of the entries.
For a set of integers $M=\{r_1, \ldots , r_m\}$ we denote by
$\S(M)$ or $\S(r_1, \ldots, r_m)$ the Smith normal form of the
diagonal matrices with the entries $r_1, \ldots , r_m$ on the diagonal.

\begin{theorem}\label{snf}
Let $p$ be a prime, and let $\X_{RC}$ be the $p$-regular character table
of $S_n$. Then we have
$$\S(\X_{RC}) = \S(b_{\mu} \mid \mu  \in C)_{p'} \: .$$
\end{theorem}

\textbf{Proof.}
Let $\Y=\Y_{CC}=(\xi^{\lambda}(\mu))_{\lambda,\mu \in C}$
denote the part of the permutation character table of $S_n$
with rows and columns indexed by the class $p$-regular partitions of~$n$.
Set $\X=\X_{RC}$.\\
As the characters $\chi^{\lambda}$ with $\lambda$ in the set $R$
of $p$-regular partitions of~$n$ form a basic set for the characters
on the $p$-regular conjugacy classes by~\cite{KOR}, we have a
decomposition matrix $D=D_{CR}$ with integer entries such that
$$\Y = D \cdot \X \:.$$
Now by Corollary~\ref{permchartable} the permutation character table
$\Y$ is (with respect to a suitable ordering) a lower triangular
matrix with the $b_{\mu}$, $\mu \in C$, on the diagonal. Hence using
\cite[Theorem~4]{OlsReg} and Theorem~\ref{det-regchartable} we obtain
$$\det (\Y)_{p'} = (b_C)_{p'} = a_C = |\det (\X)| \: .$$
Thus $\det (D)$ is a $p$-power, and hence $\det (D)$ and $\det (\X)$
are coprime. 
This implies by \cite[Theorem II.15]{Newman}
$$\S(\Y) = \S(D \, \X) = \S(D) \S(\X) \: .$$
Now using the divisibility property in Corollary~\ref{permchartable} (iii)
we can convert the triangular matrix $\Y$ by unimodular transformations
to a diagonal matrix with the same entries $b_{\mu}$, $\mu \in C$, on
the diagonal, and hence $\S(\Y)=\S(b_{\mu} \mid \mu  \in C)$. As
$\S(D)$ is a diagonal matrix with only $p$-power entries on the
diagonal, this yields the assertion in the Theorem. 
$\ \Box$

\bigskip

{\bf Remark.}
Choosing $p>n$ in Theorem~\ref{snf}
shows in particular that the Smith normal form
of the whole character table $\X$ is the same
as that of the diagonal matrix with diagonal entries $b_{\mu}=R_{\mu\mu}$,
$\mu\in \P$.
One may also use the language of
symmetric functions to prove this result.
Here, one uses that the matrix $\X$ is the transition matrix from the
Schur functions to the power sums \cite[Cor.~7.17.4]{ec2}.
 Since the transition matrix from the monomial
symmetric functions to the Schur functions is an integer matrix of
determinant 1 (in fact, lower unitriangular with respect to a
suitable ordering on partitions \cite[Cor.~7.10.6]{ec2}), the
transition matrix $R_n=(R_{\lambda\mu})_{\lambda,\mu\in \P}$
between the $m_\lambda$'s and $p_\mu$'s has the same  Smith normal form
as $\X$. Then we use the same arguments as before to deduce
the Smith normal form of $R_n$.

\bigskip

{\bf Remark.}
We do not know at present how Theorem \ref{snf} should extend from the
prime case to the case of general $r.$
Some obvious guesses for $r$-versions do not hold.
The following weaker version {\it might} be true. Let $\pi$ be the set of
primes of $r$, and for a number $m$ let $m_{\pi'}$ denote its
$\pi'$-part (the largest divisor of $m$ coprime to $r$).
Then
$$\S(\X_{RC})_{\pi'}= \S(b_{\mu} \mid \mu \in C)_{\pi'} \: .$$

\bigskip

Using Theorem~\ref{snf} above for $p=2$ also allows the determination
of the Smith normal form of the reduced spin character table of the
double covers of the symmetric groups. For the background on spin
characters of $S_n$ we refer to \cite{HoHu} and \cite{Schur}. 

We denote by $\D$ the set of partitions of $n$ into distinct parts and
by $\O$ the set of partitions of $n$ into odd parts.  Note that thus
$\D$ is the set of 2-regular partitions of $n$ and $\O$ is the set of
class 2-regular partitions of $n$.  For each $\lambda \in \D$ we have
a spin character $\gen{\lambda}$ of $S_n$.  If $n-\ell(\lambda)$ is
odd, then there is an associate spin character $\gen{\lambda}'=\sgn
\cdot \gen{\lambda}$ of $S_n$ and $\lambda$ is said to be of negative
type; the corresponding subset of $\D$ is denoted by $\D^-$.  The spin
characters can have non-zero values only on the so-called doubling
conjugacy classes of the double cover $\tilde{S}_n$ of $S_n$; these
are labelled by the partitions in $\O \cup \D^-$.  More precisely, for
any such partition we have two conjugacy classes in $\tilde{S}_n$; one
of these is chosen in accordance with \cite{Schur}, and we denote a
corresponding representative by $\sigma_{\mu}$.  While the spin
character values on the $\D^-$ classes are known explicitly (but they
are in general not integers, and mostly not even real), for the values
on the $\O$-classes we only have a recursion formula (due to A.\
Morris) which is analogous to the Murnaghan-Nakayama formula, and
which shows that these are integers.  We then define the reduced spin
character table as the integral square matrix
$$Z_s = (\gen{\lambda}(\sigma_{\mu}))_{{\lambda\in \D} \atop {\mu \in
    \O}}$$

For any integer $m\ge 0$, let $s(m)$ be the number of summands in the
2-adic decomposition of $m$.
For $\alpha=(1^{m_1},3^{m_3},\cdots)\in \O$ we define
$$k_\alpha=\sum_{i\;{\mathrm{odd}}}(m_i-s(m_i))\:.$$
Then we have

\begin{theorem}\label{spin-snf}
The Smith normal form of the reduced spin character table $Z_s$ of $\tilde{S}_n$
is given by
$$\S(Z_s)=\S(2^{[k_{\mu}/2]}, \mu \in \O) \cdot \S(b_{\mu}, \mu \in \O)_{2'} \; .$$
\end{theorem}

{\bf Proof.}
Let $\Phi$ denote the Brauer character table of $\tilde{S}_n$ at characteristic~2;
this is equal to the Brauer character table of $S_n$.
Then
$Z_s=D_s \cdot \Phi$,
where $D_s$ is a ``reduced'' decomposition matrix at $p=2$;
the reduction corresponds to leaving out the associate spin characters $\gen{\lambda}'$
for $\lambda \in \D^-$.
The matrix $D_s$ is then an integral square matrix.
In \cite{BesOls2pow}, the Smith normal form of $D_s$ was determined:
$$\S(D_s)=\S( 2^{[k_{\mu}/2]}, \mu \in \O)\: .$$
As this is a matrix of 2-power determinant and the determinant of the Brauer character table
is coprime to~2, we have
$$\S(Z_s)=S(D_s) \cdot \S(\Phi) = \S( 2^{[k_{\mu}/2]}, \mu \in \O) \cdot \S(\Phi)\: .$$
Now the Brauer characters and the characters $\chi^{\lambda}, \lambda \in R=\D$,
are both basic sets for the characters of $S_n$ on 2-regular classes, hence
$\S(\Phi)=S(\X_{RC})$.
By Theorem~\ref{snf} (for $p=2$) we thus obtain
$$\S(\Phi)=\S(\X_{RC}) = \S(b_{\mu} \mid \mu  \in \O)_{2'} \: .$$
This proves the claim.
$\ \Box$

\bigskip

{\bf Remark.}
Let us finally mention some open questions.
We have determined the Smith normal form for the whole
reduced spin character table.
It is natural to ask whether also a $p$-version (or even an $r$-version)
of this holds, or at least, whether  the determinant can be computed
similarly as in the ordinary $S_n$ case.

More precisely, for a prime $p$ define
$$Z_{s,p}=(\gen{\lambda}(\sigma_{\mu}))_{{\lambda\in \D_p} \atop {\mu \in \O_p}}$$
where $\D_p$ and $\O_p$ denote the sets of class $p$-regular partitions in $\D$
and $\O$, respectively.
Some examples lead to the following conjecture:
$$\S(Z_{s,p})=\S(2^{[k_{\mu}/2]}, \mu \in \O_p) \cdot \S(b_{\mu}, \mu \in \O_p)_{2'} \; .$$

Concerning the determinant, one may ask whether
there is an analogue of Theorem~\ref{det^2-regchartable} in the spin case.

For $S_n$ as well as its double cover one may also try to look for
 sectional versions or block versions
for the results on regular character tables.

\end{document}